
%
%
\documentstyle{amsppt}
\magnification \magstep1
\parskip 6pt
\parindent .3in
\pagewidth{5.2in}
\pageheight{7.55in}
\font\bbewend=lasy10

\def\MP #1{{\Bbb P}^#1}

\def\lra {{\longrightarrow}}
\def\bewende {{\vskip-5pt {\hfill {\bbewend \char50}}\vskip 4pt}\noindent}

\def\QQ {{\Bbb Q }}
\topmatter

\title Ample line bundles on blown up surfaces
\endtitle
\author Oliver K\"uchle \endauthor
\address Max-Planck-Instutut f\"ur Mathematik, Gottfried-Claren-Str. 26, 53225
Bonn
\endaddress \email kuechle\@mpim-bonn.mpg.de \endemail
\thanks
Supported by Max-Planck-Gesellschaft
\endthanks
\rightheadtext{} \leftheadtext{}
\abstract
Given a smooth complex projective surface $S$
and an ample
divisor $H$ on $S$, consider the blow up of $S$ along $k$
points in general position. Let $H'$ be the pullback of
$H$ and $E_1,..., E_k$ be the exceptional divisors.
We show that $L = nH' -E_1 - ... -E_k$ is ample
if and only if $L^2$ is positive  provided the integer $n$ is at least 3.
\endabstract

\endtopmatter

\document

\subhead
Introduction
\endsubhead

In this note we give an answer to the following question:
Given a smooth projective surface $S$ over ${\Bbb C}$ and an ample
divisor $H$ on $S$, consider the blow up $f: S'\lra S$ of $S$ along $k$
points in general position. Let $H'=f^*H$
and $E_1,\dots, E_k$ be the exceptional divisors.
When is the divisor
$$L= nH' -\sum_{i=1}^k E_i$$
ample ?
\par\noindent
We show that the condition $L^2 > 0$, which clearly
is necessary, is also sufficient provided the integer $n$ is at least 3.
Note that the answer to this question has been unknown even
in the case of $S=\MP2$.
The basic idea is to study the situation on the surface $S$
with variational methods.
\par\noindent
Shortly after this work has been completed the author learned that Geng
Xu obtained a similar result in the case of $S=\MP2$ independently.
\smallskip\noindent
It's a pleasure to thank Rob Lazarsfeld, who introduced me
to this circle of ideas.
\subhead
Proofs
\endsubhead

The main technical tool is an estimate on the self-intersection of moving
singular curves established by Ein, Lazarsfeld and Xu in the context of
Seshadri constants of ample line bundles on smooth surfaces (cf. \cite{EL},1.2,
and \cite{Laz}, 5.16). The precise statement is:
\proclaim{Proposition}
Let $\{ C_t \}_{t\in\Delta}$ be a 1--parameter family of reduced irreducible
curves on a smooth projective surface $X$, and $y, y_1,\dots ,y_r\in X$
be distinct points such that $mult_{y_i} C_t\ge m_i$ for all $t\in\Delta$
and $i=1,\dots ,r$. Suppose there exist $t, t'$ with
$mult_y C_{t}=m >0$ and $y\not\in C_{t'}$. Then
$$(C_t)^2\ge m(m-1) +\sum_{i=1}^r m_i^2.$$
\endproclaim\bewende \noindent
Using this Proposition we can prove:
\proclaim{Theorem}
Let $S'$ be as above and $a >2$ be a rational
number. Consider the $\QQ$-divisor
$$M= aH' -\sum_{i=1}^k E_i.$$
Then the following hold:
\roster
\item
If  $M^2= a^2 H^2 - k\ge 2$, then $M$ is ample on $S'$.
\medskip
\item
If $M^2 =a^2 H^2 -k\ge 1$, then $M$ is positive on all curves
$C'\subset S'$ for which $j$ exists with $C'.E_j\ge 2$.
\endroster
\endproclaim
\demo{Proof}
Suppose the theorem is not true, and choose
 a curve $C'\subset S'$ such that $M.C' \le 0$.
Consider $C =f(C')$. Defining $m_i=mult_{p_i}(C)$, we may suppose
that $m_1\ge\dots\ge m_k$.
Since $M.C' \le 0$, we have
$$\sum_{i=1}^k m_i \ge a(H.C).\eqno{(*)}$$
Now we may assume that
\roster
\item"--"
$C$ passes through all the points $p_i$, i.e. $m_i\ge 1$\smallskip
\item"--"
$C$ is irreducible and reduced\smallskip
\item"--"
$C$ moves, since the $p_i$ are in general position\smallskip
\endroster
Here $C$ moves even in the strong sense, that is, fixing $p_1,\dots ,p_{k-1}$,
the curve
$C$ still moves in a family of curves satisfying $(*)$.
To see this simply observe that
any curve on $S$ lies in one of countably many families,
but no neighbourhood of $p_k$ is covered by countably many curves.\par\noindent
Finally we claim that a general member of this family has sufficiently
big multiplicity at $p_1,\dots , p_{k-1}$. But  any member
satisfies $(*)$, so this follows from semicontinuity.
\medskip\noindent
Therefore  we can apply the Proposition and
obtain the estimate
$$C\cdot C\ge m_1^2+\dots + m_{k-1}^2+ m_k(m_k-1),$$ and hence
combined with  the Hodge-Index-Theorem
$$\Biggl(\sum_{i=1}^k m_i\Biggr)^2 \ge a^2(H.C)^2\ge a^2 H^2\cdot C^2\ge
 a^2 H^2\Bigl(\sum_{i=1}^k m_i^2 -m_k\Bigr).\eqno{(**)}$$
By $(*)$, $(**)$ and the assumption $a > 2$ we may assume
$k\ge 2$ in the following.\par\noindent
Suppose for the time being that $C$ is not smooth at one of the $p_j$,
which is the case if and only if $C'.E_j\ge 2$. Then $m_1\ge 2$, and
$(**)$
contradicts the following Lemma:
\enddemo
\proclaim{Lemma}
Let $k\ge 2$ and $x_1,\dots ,x_k\in {\Bbb Z}$ be integers with
$x_1\ge\dots\ge x_k\ge 1$ and $x_1\ge 2$.
Then we have
$$(k+1)\sum_{i=1}^k x_i^2 > \Bigl(\sum_{i=1}^k x_i\Bigr)^2+x_k(k+1).$$
\endproclaim
\demo{Proof of the Lemma}
We argue by induction on $k\ge 2$.
\par\noindent
For $k=2$ one proves
$$3(x_1^2+x_2^2)- (x_1+x_2)^2 - 3 x_2 > 0$$
by minimizing this expression with respect to $x_2$.
{}From the inductive hypothesis, we then obtain
$$(k+1)\sum_{i=1}^k x_i^2 >
kx_k^2+\sum_{i=1}^k x_i^2+\Biggl(\sum_{i=1}^{k-1} x_i\Biggr)^2+k x_k$$
$$ =\Biggl(\sum_{i=1}^k
x_i\Biggr)^2+x_k(k+1)-x_k^2-2\cdot\sum_{i=1}^{k-1}x_ix_k -
x_k +k x_k^2+\sum_{i=1}^k x_i^2$$
$$=\Biggl(\sum_{i=1}^k x_i\Biggr)^2 +x_k(k+1)+
\sum_{i=1}^{k-1}(x_i-x_k)^2 + x_k^2 - x_k.$$
So what we need to show is
$$\sum_{i=1}^{k-1}(x_i-x_k)^2 + x_k^2 \ge x_k,$$
but this is obvious.
\enddemo
\bewende
This proves the second
part of the Theorem.
To prove the first part it remains to exclude the case $m_1=\dots =m_k=1$.
But then $(**)$ reads
$$k^2 \ge H^2\cdot a^2(k-1),$$ contradicting the assumptions on $a$.
\bewende
\medskip\noindent
\proclaim{Corollary}
Let $L$ be as in the introduction. Then $L$ is ample if and only if $L^2 > 0$.
\endproclaim
\demo{Proof}
It clearly suffices to prove the if--part.
So  suppose $L^2 >0$ and that $L$ is not ample. Then by the Theorem
we know that $L^2=1$, i.e. $n^2 H^2= k+1$, and that there
exists an irreducible reduced curve
$C\subset S$ which is smooth at all the $p_i$ satisfying
$k\ge n (H.C)$.
\par\noindent
We claim that $k=n(H.C)$ holds. Otherwise we have $L.C' <0$.
Consider the surface $\Hat S$
obtained from $S'$ by contracting the exceptional divisor $E_j$, where
$j$ is an index such that $C$ passes smoothly through $p_j$.
The image $\Hat L$ of $L$ then satisfies $\Hat L^2=L^2+1 = 2$,
hence it is ample by the Theorem.
But this contradicts $L.C' +1= \Hat L.\Hat C \le  0$
for the image $\Hat C$ of $C'$.
\smallskip\noindent
Therefore we conclude $k+1 = n^2 H^2 = n(H.C) +1$, but this is
impossible since besides $n\ne 1$ also $H^2$ and $(H.C)$ are integers.
\bewende
\enddemo
\subhead Remark  \endsubhead
The example of a line in $\MP2$ through any two points
shows that we cannot drop the assumption $n\ge 3$ in general.
On the other hand an analysis of the proof shows that
the Corollary still holds in the case $n\ge 2$ if two general
points on $S$ can not be joined by a curve $C$ with $(H.C)=1 $,
which is true e.g. whenever $H^2\ge 2$.
 \bigskip
\head
References
\endhead

\enddocument